\documentclass[12pt,leqno,fleqn,epsfig]{article}
\usepackage{amssymb,epsfig, amsmath, amsthm}
\usepackage{mathrsfs}       

\newcommand{\R}{\mathbb{R}}

\newcommand{\N}{\mathbb{N}}

\newcommand{\ba}{\boldsymbol a}
\newcommand{\bb}{\boldsymbol b}
\newcommand{\bc}{\boldsymbol c}
\newcommand{\bd}{\boldsymbol d}
\newcommand{\bfe}{\boldsymbol e}
\newcommand{\bbf}{\boldsymbol f}
\newcommand{\bg}{\boldsymbol g}
\newcommand{\bh}{\boldsymbol h}
\newcommand{\bi}{\boldsymbol i}
\newcommand{\bj}{\boldsymbol j}
\newcommand{\bk}{\boldsymbol k}
\newcommand{\bl}{\boldsymbol l}
\newcommand{\bm}{\boldsymbol m}
\newcommand{\bn}{\boldsymbol n}
\newcommand{\bo}{\boldsymbol o}
\newcommand{\bp}{\boldsymbol p}
\newcommand{\bq}{\boldsymbol q}
\newcommand{\br}{\boldsymbol r}
\newcommand{\bs}{\boldsymbol s}
\newcommand{\bt}{\boldsymbol t}
\newcommand{\bu}{\boldsymbol u}
\newcommand{\bv}{\boldsymbol v}
\newcommand{\bw}{\boldsymbol w} 
\newcommand{\bx}{\boldsymbol x}
\newcommand{\by}{\boldsymbol y}
\newcommand{\bz}{\boldsymbol z}
\newcommand{\buu}{ \underline{\boldsymbol u}}

\newcommand{\bA}{\boldsymbol A}
\newcommand{\bB}{\boldsymbol B}
\newcommand{\bC}{\boldsymbol C}
\newcommand{\bD}{\boldsymbol D}
\newcommand{\bE}{\boldsymbol E}
\newcommand{\bF}{\boldsymbol F}
\newcommand{\bG}{\boldsymbol G}
\newcommand{\bH}{\boldsymbol H}
\newcommand{\bI}{\boldsymbol I}
\newcommand{\bJ}{\boldsymbol J}
\newcommand{\bK}{\boldsymbol K}
\newcommand{\bL}{\boldsymbol L}
\newcommand{\bM}{\boldsymbol M}
\newcommand{\bO}{\boldsymbol O}
\newcommand{\bP}{\boldsymbol P}
\newcommand{\bQ}{\boldsymbol Q}
\newcommand{\bR}{\boldsymbol R}
\newcommand{\bS}{\boldsymbol S}
\newcommand{\bT}{\boldsymbol T}
\newcommand{\bU}{\boldsymbol U}
\newcommand{\bV}{\boldsymbol V}
\newcommand{\bW}{\boldsymbol W}
\newcommand{\bX}{\boldsymbol X}
\newcommand{\bY}{\boldsymbol Y}
\newcommand{\bZ}{\boldsymbol Z}

\newcommand{\bfvarphi}{\boldsymbol\varphi}
\newcommand{\bfeta}{\boldsymbol\eta}

\newcommand{\bfpsi}{\boldsymbol\psi}

\newcommand{\bLambda}{\boldsymbol \Lambda}

\newcommand{\mes}{\operatorname{\rm meas}}

\newcommand{\supp}{\operatorname*{supp}}

\newcommand{\const}{\operatorname*{const}}

\newcommand{\curl}{\operatorname{curl}}

\newcommand{\be}{\begin{equation}}
\newcommand{\ee}{\end{equation}}
\newcommand{\bea}{\begin{eqnarray}}
\newcommand{\eea}{\end{eqnarray}}
\newcommand{\bean}{\begin{eqnarray*}}
\newcommand{\eean}{\end{eqnarray*}}
\newcommand{\supl}{\sup\limits}
\newcommand{\var}{\varepsilon}

\newcommand{\intl}{\int\limits}

\newcommand{\Beweisende}{\rule{0.2cm}{0.2cm}}

\newcommand{\D}{\displaystyle}

\newcommand{\intmw}{{\int\hspace{-830000sp}-\!\!}}

\newcounter{secnum}

\newtheorem{thm}{Theorem}[section]

\newtheorem{lem}[thm]{Lemma}

\theoremstyle{definition}
\newtheorem{defin}[thm]{Definition}
\newtheorem{rem}[thm]{Remark}

\title{Regularity  of the   $3D$ stationary  Hall magnetohydrodynamic equations
on the plane
}
 \author{Dongho Chae$^*$  and J\"{o}rg Wolf $^\dagger$\\
\ \\
 $*$Department of Mathematics\\
Chung-Ang University\\
 Seoul 156-756, Republic of Korea\\
 e-mail: dchae@cau.ac.kr\\
and \\
$\dagger$Department of Mathematics\\
 Humboldt University of Berlin\\                                           
Unter den Linden 6, 10099 Berlin, Germany\\
e-mail: jwolf@math.hu-berlin.de}
\date{}
\begin{document}
\maketitle
\begin{abstract}
We study the  regularity of weak solutions to  the  3D valued stationary Hall magnetohydrodynamic equations on $ \Bbb R^2$. We prove that every weak solution  is smooth. Furthermore,  we prove a Liouville type theorem for the Hall equations.  
  \noindent{\bf
keywords:}  stationary Hall-MHD equations, regularity, Liouville type theorem

\end{abstract}


\section{Introduction and the main theorems}
\label{sec:-1}
\setcounter{secnum}{\value{section} \setcounter{equation}{0}
\renewcommand{\theequation}{\mbox{\arabic{secnum}.\arabic{equation}}}}

We study  the following $3D$ valued stationry Hall-magnetoydrodynamics(Hall-MHD) system on $\R^2$.
\begin{align}
 (\bv\cdot \nabla )\bv - \Delta \bv &= -\nabla p +
 (\nabla \times \bB) \times  \bB + \bbf,
\label{1.1}
\\
 \nabla \times (\bB\times \bv) - \Delta \bB &= 
-\nabla \times ((\nabla \times \bB) \times  \bB)+  \nabla \times \bg,
\label{1.2}
\\
\nabla \cdot \bv &=0,\quad  \nabla \cdot \bB =0.
\label{1.3}  
  \end{align}
Here, $\bv=(v^1, v^2, v^3), \bB=(B^1, B^2, B^3)$, where $v^j=v^j (x_1, x_2), B^j=B^j(x_1,x_2), j=1,2,3$, and $p=p(x_1,x_2)$, $x=(x_1,x_2)\in \R^{2}$.  The vector fields $\bbf$ and $\bg$ represent the external forces. The system for $\bB$ obtained from (\ref{1.2}) by setting $\bv=0$ is called the Hall equation.
Physically the full time-dependet version of the system (\ref{1.1})-(\ref{1.3}) on $\Bbb R^3$ describes the dynamics plasma flows with strong shear of magnetic fields such case as in the solar flares. We refer \cite{ach} and the references therein for the physical backgrounds for the full system, and  \cite{cha1, cha2, cha3, cha4, cha5, cha6, sue, fan} for  recent studies of the  mathematical problems of the equations.  In particular in \cite{cha6} it is shown that there exist weak solutions of the time dependent 3D  time dependent Hall-MHD system on the plane, having  the  possible set of space-time singularities, whose Hausdorff dimension is at most two.   On the other hand, in  \cite{cha5} it is proved that there exist  weak solutions of  the full 3D stationary Hall-MHD equations having  the  possible set of singularities with the  Hausdorff dimension at most one. 
In the case of our system (\ref{1.1})-(\ref{1.3}) of the 3D stationary Hall-MHD on the plane, if we  apply the argument of \cite{cha5}, then we could easily deduce that there exist  weak solutions having  possible set of singularities with the Hausdorff dimension zero. Note that this is still far from the conclusion that the set of singularities is empty. Therefore, the regularity problem of (\ref{1.1})-(\ref{1.3}) in $\Bbb R^2$  could be regarded as an interesting {\em critical problem}, which is our main subject of study in this paper. 
One of our main results in this paper is to show the {\em full regularity} of any weak solutions to the above system, namely the set of singularities is indeed empty. For the proof of this result we  modify the Widman's hole filling method (cf. \cite{wid}) in order  to handle the case,  where the logarithmically blowing-up coefficient is allowed in the Caccioppoli type inequality. We also prove a Liouville type result for the Hall system, which means  that  any weak solution the equations (\ref{1.2})-(\ref{1.3}) with  $\bv=\nabla \times \bg=0$ having finite Dirichlet integral is zero.

\hspace{0.5cm}
Below by $ \nabla ^{ \bot}$ we denote the orthogonal gradient operator  $(-\partial _2, \partial _{ 1})^{ \top}$. We also denote that $ \bv' = (v^1, v^2) $ for given $ \bv = (v^1, v^2, v^3)$. 
According to 
$ \nabla \cdot \bB =0$ we 
find a potential $ \Phi\in W^{2,\, 2}_{\rm loc}(\R^{2})$ such that 
\begin{equation}
(B^1, B^2)^{ \top} = \nabla ^{ \bot} \Phi.
\label{1.p}
\end{equation}
Setting $ \Psi := B_3$,  the equations in  \eqref{1.2} turn into  
\begin{align}
\Delta \Phi &= \nabla ^{ \bot} \Psi \cdot \nabla \Phi+ h^3, 
\label{1.2a}
\\
\Delta  \Psi &= - \nabla^{ \bot} \Delta \Phi\cdot  \nabla \Phi + \partial _1 h^2 -\partial _2 h^1,
\label{1.2b}
\end{align}
where 
\begin{align}
\begin{cases}
h^1 = -\partial _1\Phi  v^3 +\Psi  v^2 + g^1,
\\
h^2 = -\partial _2\Phi  v^3 - \Psi  v^1 + g^2,
\\
h^3 =\nabla \Phi \cdot \bv'  + g^3 + \const.
\end{cases}
\label{1.2c}
\end{align}
We call the system \eqref{1.2a} - \eqref{1.2b} the $ \Phi $-$ \Psi $-{\it system}.   

\vspace{0.5cm}  
\hspace{0.5cm}
By $ {\hat \bW} ^{1,\, 2}_{ \sigma }(\R^{2})$ we denote the space of all $ \bv \in \bL^{2}_{ \rm loc}(\R^{2})$
with $ \nabla \bv \in \bL^2(\R^{2})$ and $ \nabla \cdot \bv =0$ almost everywhere in $ \R^{2}$.  In addition, by  
$ {\hat W} ^{m,\, s}(\R^{2})$ we denote the space of all $ \Phi  \in W^{m,\, s}_{ \rm loc}(\R^{2})$ with $ D^\alpha  \Phi 
\in L^2(\R^{2})$ for all $| \alpha | = m$. 

\hspace{0.5cm}
We introduce the following 
notion of   weak solution to \eqref{1.1}--\eqref{1.3}, and the notion of  weak-strong solution to the system \eqref{1.2b}, \eqref{1.2a}.  

\begin{defin}
1. Let $ \bbf , \bg \in \bL^2(\R^{2})$. A pair    $ (\bv, \bB) \in {\hat \bW}_\sigma  ^{1,\, 2}(\R^{2})\times {\hat \bW}_\sigma^{1,\, 2}(\R^{2})$ 
is called a {\it weak solution to} \eqref{1.1}--\eqref{1.3} if the following identities  hold for 
all  $ \bfvarphi \in \bC^{\infty}_{\rm c, \sigma }(\R^{2})$, and $ \bfpsi \in  \bC^{\infty}_{\rm c}(\R^{2})$ respectively 
\begin{align}
&  \intl_{\R^{2}} \nabla \bv : \nabla \bfvarphi =  \intl_{ \R^{2}}(\bv \otimes \bv ): \nabla \bfvarphi +  (\nabla \times \bB )\times \bB \cdot \bfvarphi + \bbf \cdot \bfvarphi, 
\label{1.6}
\\
 &  \intl_{\R^{2}} \nabla \bB : \nabla \bfpsi =   -\intl_{\R^{2}}  ((\nabla \times \bB) \times \bB + \bB \times \bv - \bg )\cdot  \nabla \times \bfpsi. 
\label{1.7}
\end{align}

2. Let $  \bh\in \bL ^2(\R^{2})$. A pair $ (\Phi , \Psi )\in {\hat W} ^{2,\, 2}(\R^{2})\times {\hat W} ^{1,\, 2}(\R^{2})$ 
is called a {\it strong-weak solution to} \eqref{1.2a}, \eqref{1.2b} if \eqref{1.2a} is satisfied almost everywhere in $ \R^{2}$, and  \eqref{1.2b} 
is fulfilled in the sense of distributions, i.\,e.   
for every $ \varphi \in 
C^{\infty}_{\rm c} (\R^{2})$,

\begin{align}
\intl_{\R^{2}} \nabla \Psi \cdot \nabla \varphi = \intl_{\R^{2}} -\Delta  \Phi \nabla \Phi \cdot \nabla ^{ \bot} \varphi 
 +\bh' \cdot  \nabla^{ \bot}  \varphi.
\label{1.9}
\end{align}   
 
\end{defin}

\begin{rem}
Note that if $ (\bv ,\bB) \in {\hat \bW} ^{1,\, 2}_\sigma (\R^{2})\times {\hat \bW} ^{1,\, 2}_\sigma (\R^{2})$ 
is a weak solution to \eqref{1.1}-- \eqref{1.3}, then $ (\Phi , \Psi )\in {\hat W} ^{2,\, 2}(\R^{2})\times {\hat W} ^{1,\, 2}(\R^{2})$ is a weak-strong solution to \eqref{1.2a}, \eqref{1.2b} with right-hand side $\bh$ given according to  \eqref{1.2c}.     
Indeed, noting 
\[
(\nabla \times \bB) \times \bB = - \frac{1}{2} (\nabla \Psi ^2,  0)^{ \top} - 
(\Delta \Phi \nabla \Phi , \nabla ^{ \bot} \Psi \cdot \nabla \Phi ),
\]
from \eqref{1.7} with $ \bfpsi =(\eta ^1, \eta ^2, 0)^{ \top} \in \bC^{\infty}_{\rm c}(\R^{2})$ 
we find 
\begin{align*}
 \intl_{ \R^{2}}\Delta  \Phi \curl \bfeta   &= 
\intl_{ \R^{2}} \partial _i \nabla^{ \bot} \Phi \cdot \partial _i  \bfeta 
\\
&= \intl_{\R^{2}}   (\nabla ^{ \bot} \Psi \cdot \nabla \Phi - (\bB \times \bv)^3 + g^3 ) \curl \bfeta,
\end{align*}
where $ \curl \bfeta = \partial _1 \eta ^2 - \partial _2 \eta^1$. 
Whence, \eqref{1.2a}. 

\hspace{0.5cm}
To verify \eqref{1.2b},  we insert into \eqref{1.7} the test functions $ \bfpsi =(0,0, \varphi )^{ \top}$,  $\varphi \in C^{\infty}_{\rm c}(\R^{2})$. This gives 
\begin{align*}
 \intl_{\R^{2}} \nabla \Psi \cdot \nabla \varphi  = \intl_{\R^{2}} \Delta \Phi \nabla \Phi \cdot \nabla^{ \bot} \varphi   
 +(\bB \times \bv  - \bg )' \cdot \nabla^{ \bot} \varphi,
\end{align*}
and therefore  \eqref{1.9} holds with $ \bh$ given by  \eqref{1.2c}.  Accordingly, the pair $ (\Phi , \Psi )$ is a strong-weak solution to \eqref{1.2a}, \eqref{1.2b}.   

\end{rem}

\hspace{0.5cm}
Our first main result is the following regularity theorem for the system (\ref{1.1})-(\ref{1.3}).

\begin{thm}
\label{thm1.2}
Let $ (\bv, \bB) \in {\hat \bW}_\sigma  ^{1,\, 2}(\R^{2})\times {\hat \bW}_\sigma^{1,\, 2}(\R^{2})$ be a weak solution to the 
steady Hall-MHD system in $ \R^{2}$ with $\bbf, \bg \in C^\infty (\Bbb R^2)$. Then both $ \bv $ and $ \bB $ are smooth.  
\end{thm}

\hspace{0.5cm}
Next, we  consider the following  stationary Hall system,
\begin{align}
\Delta \bB = \nabla \times ((\nabla \times \bB) \times \bB )\quad  \text{ in}\quad  \R^{2},
\label{2.0}
\end{align}
which is obtained from the $\bB$ equations of the Hall-MHD system with $ \bv \equiv {\bf 0}$.  
Our second main result is the following Liouville type theorem for the system (\ref{2.0}).

\begin{thm}
\label{thm1.4}
Let $ \bB$ be  a weak solution to \eqref{2.0} having the finite Dirichlet integral, i.e. 
$ \intl_{ \R^{2}} | \nabla \bB |^2 < +\infty.$   
Then $ \bB \equiv {\bf 0}$. 
\end{thm}

\section{A modified hole filling method}
\label{sec:-7}
\setcounter{secnum}{\value{section} \setcounter{equation}{0}
\renewcommand{\theequation}{\mbox{\arabic{secnum}.\arabic{equation}}}}

\begin{thm}
\label{thm7.1}
Let $ f\in {\hat W} ^{1,\, 2}(\R^{2})$, and let $ \mu \in (0,1)$. Suppose that for  all $ B_r \subset \R^{2}$, 
$ 0< r< \frac{1}{2}$ the following inequality holds true 
\begin{equation}
\intl_{B_{ r/2}} | \nabla f|^2 \le c_0  (1+ | (f)_{ B_r}|) \intl_{B_r  \setminus B_{ r/2}} | \nabla f|^2  + c_1 r^{ \mu }, 
\label{7.1}
\end{equation}
where $ c_0, c_1$ are positive  constants. Then $ f$ is H\"older continuous.  

\end{thm}

{\bf Proof}: 1. In view of \cite[Lemma\,2.4]{fre} we see that for all $ 0< r< \frac{1}{2}$, 
\begin{equation}
1+  | (f)_{ B_r } |  \lesssim  \left(\log r^{ -1}\right)^{\frac12}.
\label{7.2}
\end{equation}

According to \eqref{7.1} together with \eqref{7.2} we find that there exists  a constant $ c>0$ such that for all $ 0< r< \frac{1}{2}$. 
\begin{equation}
\intl_{B_{ r/2}} | \nabla f |^2   \le c \left(\log r^{ -1}\right)^{\frac12}
\intl_{B_r  \setminus B_{ r/2 }} | \nabla f |^2 + c r^\mu.  
\label{7.3}
\end{equation}
Now in \eqref{7.3} filling  the hole by adding $ c \left(\log r^{ -1}\right)^{\frac12} \intl_{B_{ r/2}} | \nabla f |^2 $ to both sides, 
we infer  
\[
\intl_{B_{ r/2}} | \nabla f |^2   \le \frac{c \left(\log r^{ -1}\right)^{\frac12}}{1+c \left(\log r^{ -1}\right)^{\frac12} }
\intl_{B_r  } | \nabla f |^2  + c r^\mu.  
\]
Accordingly, we are in a position to apply  Lemma\,\ref{lemA.7}, which yields for all $ 0<r< \frac{1}{2}$ and $\alpha >1$,
\begin{equation}
\intl_{B_r} | \nabla f|^2  \lesssim  \frac{1}{[\log r^{ -1}]^{ 2\alpha } }.   
\label{7.5}
\end{equation}

2. Next, applying Lemma\,\ref{lemA.5} , we conclude that 
\[
\sup_{ 0<r<1} | (f)_{ B_r}|  \lesssim  \zeta (\alpha)<+\infty,  
\]
where the hidden constant in this inequality is independent of the center of the ball. 
Thus observing \eqref{7.1}, we get a constant $ c_2>0$  such that for all $ 0< r< \frac{1}{2}$,
\begin{equation}
\intl_{B_{ r/2}} | \nabla f|^2 \le c_2 \intl_{B_r  \setminus B_{ r/2}} | \nabla f|^2  + c_1 r^{ \mu }.
\label{7.9}
\end{equation}
Now in \eqref{7.9} filling the hole, we arrive at 
\begin{equation}
\intl_{B_{ r/2}} | \nabla f|^2 \le \theta  \intl_{B_r } | \nabla f|^2  + c_1 r^{ \mu }, \quad  
\text{ where}\quad  \theta =\frac{c_2 }{1+ c_2 }<1. 
\label{7.10}
\end{equation}

3. Let  $ 0< \lambda < \min \Big\{- \frac{\log \theta }{\log2}, \mu \Big\}$ arbitrarily chosen but fixed. 
Thanks to Lemma\,\ref{lemA.8} we get constant $ c_3 >0$ such that for all $ 0< r< \frac{1}{2}$
\begin{equation}
\intl_{B_{ r/2}} | \nabla f|^2 \le c_3 r^{ \lambda }. 
\label{7.11}
\end{equation}
Note that $ c_3$ depends neither on $ r$ nor on the center of the ball. 

\vspace{0.2cm}
4. Finally, applying Poincar\'e's inequality from \eqref{7.11} we conclude that for all $ 0< r< \frac{1}{2}$
\[
\intl_{B_r} | f- (f)_{ B_r}| ^2  \lesssim  r^2\intl_{B_{ r/2}} | \nabla f|^2 \le c_3 r^{2+ \lambda }. 
\] 
By Campanato's theorem (see e.g. \cite{gia}) we get the H\"older continuity of $ f$.  \hfill \Beweisende 

\section{Local energy equality for  weak solutions to the $ \Phi$- $\Psi$-system}
\label{sec:-8}
\setcounter{secnum}{\value{section} \setcounter{equation}{0}
\renewcommand{\theequation}{\mbox{\arabic{secnum}.\arabic{equation}}}}

The aim of this section is to show that every weak-strong solution to \eqref{1.2a},\eqref{1.2b} satisfies a corresponding local energy equality.  
We have the following 

\begin{lem}
\label{lem8.1}
Let $ \bh \in \bL ^2(\R^{2})$. Let $(\Phi , \Psi ) \in {\hat W} ^{2,\, 2}(\R^{2})\times {\hat W} ^{1,\, 2}(\R^{2}) $
be a strong-weak solution   to the $ \Phi $-$ \Psi $ system 
\eqref{1.2a}, \eqref{1.2b}.  Then the following energy identity holds true for all $ \zeta \in C^{\infty}_{\rm c}(\R^{2})$, and 
for all $ c\in \R$ 
\begin{align}
& \intl_{ \R^{2}} ((\Delta \Phi )^2 + | \nabla \Psi |^2) \zeta 
\cr
&=  -\intl_{ \R^{2}}  \Big((\Psi - c) \nabla \Psi - (\Psi - c)\Delta \Phi \nabla ^{ \bot} \Phi \Big) \nabla \zeta 
+ \intl_{ \R^{2}}  h^3 \Delta \Phi \zeta  + \bh'\cdot 
\nabla^{ \bot} ((\Psi -c)\zeta).  
\label{8.1}
\end{align}

\end{lem}

{\bf Proof}: For  $ \rho  >0$ we define 
\[
\gamma_\rho  (\tau ) = \begin{cases}
\frac{1}{2}\tau^2 \quad  &\text{ if} \quad  | \tau | \le  \rho  
\\[0.3cm]
\rho \Big(| \tau| - \frac{r}{2}\Big) \quad   &\text{ if} \quad  | \tau | >  \rho .    
\end{cases}
\]
Clearly, $ \gamma \in C^{ 1,1}(\R)$, and   
\[
\gamma '_\rho  (\tau )   = {\rm sign}(\tau )\min \{ | \tau |, \rho\},
\quad   \gamma_\rho  ''(\tau ) = \chi_{ (-\rho ,\rho )}.  
\]

\hspace{0.5cm}
Let $ \zeta \in C^{\infty}_{\rm c}(\R^{2})$ be arbitrarily chosen. By virtue of Sobolev's embedding theorem we see that 
$ \Phi$ is H\"older continuous, and thus bounded on $ \supp(\zeta )$. Without loss of generality we may assume that 
$ \Phi \ge 1$   on $ \supp(\zeta )$. Let $ \alpha >0$.  We multiply \eqref{1.2a} by $ \alpha \Phi^{ \alpha -1} \gamma_\rho  (\Psi )\zeta $, integrate  it over $ \R^{2} $,  and integrate by part. This leads to the following identity 
\begin{align}
\alpha\intl_{\R^{2}} \Delta  \Phi \Phi^{ \alpha -1} \gamma_\rho  (\Psi )\zeta &=  
\alpha \intl_{\R^{2}} (\nabla ^{ \bot} \Psi \cdot \nabla \Phi + h^3)\Phi^{ \alpha -1} \gamma_\rho  (\Psi )\zeta
\cr
&= -\intl_{\R^{2}}  \Phi ^{ \alpha }\gamma_\rho  (\Psi ) \nabla ^{ \bot} \Psi \cdot  \nabla \zeta 
+
\alpha \intl_{\R^{2}}\Phi^{ \alpha -1} \gamma_\rho  (\Psi ) h^3\zeta.
\label{8.2}
\end{align}  

On the other hand, applying integration by parts, we find 
\begin{align}
&\alpha\intl_{\R^{2}} \Delta  \Phi \Phi^{ \alpha -1} \gamma_\rho  (\Psi )\zeta
\cr
&\qquad = - \alpha(\alpha -1)\intl_{\R^{2}} | \nabla \Phi|^2 \Phi^{ \alpha -2} \gamma_\rho  (\Psi )\zeta
- \intl_{\R^{2}} \nabla  \Phi^{ \alpha }\cdot   \nabla \Psi \gamma'_\rho  (\Psi )\zeta
\cr
&\qquad \qquad \qquad - \alpha\intl_{\R^{2}} \Phi^{ \alpha -1} \gamma_\rho  (\Psi )  \nabla  \Phi\cdot \nabla \zeta
\cr
&\qquad = - \alpha(\alpha -1)\intl_{\R^{2}} | \nabla \Phi|^2 \Phi^{ \alpha -2} \gamma_\rho  (\Psi )\zeta
-\intl_{\R^{2}} \nabla  (\Phi^{ \alpha }\gamma'_\rho  (\Psi )\zeta)  \cdot \nabla \Psi 
\cr
&\qquad \qquad +
\intl_{\R^{2}} \Phi ^\alpha | \nabla \Psi |^2  \gamma_\rho  ''(\psi ) \zeta+ 
\intl_{\R^{2}}  \Phi ^{ \alpha }  \gamma_\rho'(\Psi ) \nabla \Psi \cdot \nabla \zeta
\cr
&\qquad \qquad \qquad - \alpha\intl_{\R^{2}} \Phi^{ \alpha -1} \gamma_\rho  (\Psi )  \nabla  \Phi\cdot \nabla \zeta.
\label{8.3}
\end{align}
In what follows, we focus  on evaluating the second integral on the right-hand side. For this purpose we  first replace $ \gamma_\rho  '(\Psi )$ by 
$ \gamma_\rho  '(\Psi )_{ \var } =\gamma_\rho  '(\Psi ) \ast \eta _\var $, where $ \eta _\var $ denotes the usual 
Friedrich's mollifying kernel. From \eqref{1.9} with $ \varphi =  \Phi ^\alpha \gamma_\rho  '(\Psi )_{ \var }$  we get 
\begin{align}
&-\intl_{\R^{2}}\nabla  (\Phi^{ \alpha }\gamma'_\rho  (\Psi )_{ \var }\zeta)  \cdot \nabla \Psi 
\cr
&\qquad = 
\intl_{\R^{2}} \Delta  \Phi \nabla \Phi \cdot \nabla ^{ \bot}(\Phi^{ \alpha }\gamma'_\rho  (\Psi )_{ \var }\zeta)
-\intl_{\R^{2}}\bh' \cdot  \nabla^{ \bot} (\Phi^{ \alpha }\gamma'_\rho  (\Psi )_{ \var }\zeta) =I_\var +II_\var . 
\label{8.4}
\end{align}   
By an elementary calculus we get 
\begin{align}
I_\var & = \intl_{\R^{2}}\Phi ^{ \alpha }\gamma_\rho ' (\Psi )_{ \var } \Delta  \Phi \nabla \Phi \cdot  \nabla ^{ \bot}\zeta
+ \intl_{\R^{2}} \Phi ^\alpha \Delta  \Phi (\nabla \Phi \cdot \nabla ^{ \bot} \Psi  \gamma''_\rho (\Psi ))_{ \var }\zeta 
\cr
& \qquad -\intl_{\R^{2}} \Phi ^\alpha \Delta  \Phi 
\Big[(\nabla \Phi \cdot \nabla ^{ \bot}\gamma'_\rho (\Psi ))_{ \var }
- \nabla \Phi \cdot (\nabla ^{ \bot}  \gamma'_\rho (\Psi )   )_{ \var }\Big] \zeta.
\label{8.4a}  
\end{align}
As it can be checked easily,  the first integral on the right-hand side of \eqref{8.4a} tends to 
\[
\intl_{\R^{2}}\Phi ^{ \alpha }\gamma_\rho ' (\Psi)  \Delta  \Phi \nabla \Phi \cdot  \nabla ^{ \bot}\zeta \quad  \text{ as}\quad 
\var \rightarrow 0,
\]
while the second integral tends to 
\[
 \intl_{\R^{2}} \Phi ^\alpha \Delta  \Phi \nabla \Phi \cdot \nabla ^{ \bot} \Psi  \gamma''_\rho (\Psi ))\zeta 
\quad  \text{ as}\quad \var \rightarrow 0,
\]
where we have used the fact that $\nabla \Phi \cdot \nabla ^{ \bot} \Psi \gamma_\rho  ''(\Psi )= \Delta \Phi \gamma_\rho  ''(\Psi )\in  L^2(\R^{2}) $. Finally, appealing to Lemma\,\ref{lem8.2} below with $ \psi = \gamma_\rho  '(\Psi )$, and $ \phi = \Phi $, 
we infer that the  third integral  tends to zero as $ \var \rightarrow 0$. 

\hspace{0.5cm}
Furthermore, the convergence of the integral $ II_\var $, and the convergence of the integral on the left-hand side of \eqref{8.4}  can be obtaind by using routine arguments, recalling  the fact that $ f_\var \rightarrow f$ 
in $ L^1(\R^{2})$ as $ \var \rightarrow 0$ for any $ L^1$ function $ f$.   This together with \eqref{1.2a} shows that 
\begin{align*}
&-\intl_{\R^{2}}\nabla  (\Phi^{ \alpha }\gamma'_\rho  (\Psi )\zeta)  \cdot \nabla \Psi 
\\
&\qquad = -
\lim_{\var  \to 0} \intl_{\R^{2}}\nabla  (\Phi^{ \alpha }\gamma'_\rho  (\Psi )_{ \var }\zeta)  \cdot \nabla \Psi 
\\
&\qquad = \intl_{\R^{2}}\Phi ^{ \alpha }\gamma_\rho ' (\Psi)  \Delta  \Phi \nabla \Phi \cdot  \nabla ^{ \bot}\zeta
+ \intl_{\R^{2}} \Phi ^\alpha ((\Delta  \Phi)^2 - h^3\Delta \Phi  ) \gamma''_\rho (\Psi ))\zeta
 \\
&\qquad \qquad   -\intl_{\R^{2}}\bh' \cdot  \nabla^{ \bot} (\Phi^{ \alpha }\gamma'_\rho  (\Psi )\zeta).  
\end{align*} 
Replacing the second integral on the right-hand side of \eqref{8.3} by the identity,  we have just derived, we obtain 
\begin{align}
&\alpha\intl_{\R^{2}} \Delta  \Phi \Phi^{ \alpha -1} \gamma (\Psi )\zeta
\cr
&\qquad = - \alpha(\alpha -1)\intl_{\R^{2}} | \nabla \Phi|^2 \Phi^{ \alpha -2} \gamma (\Psi )\zeta +
\intl_{\R^{2}}\Phi ^{ \alpha }\gamma_\rho ' (\Psi)  \Delta  \Phi \nabla \Phi \cdot  \nabla ^{ \bot}\zeta
\cr
&\qquad \qquad + \intl_{\R^{2}} \Phi ^\alpha ((\Delta  \Phi)^2 - h^3\Delta \Phi  ) \gamma''_\rho (\Psi ))\zeta
-\intl_{\R^{2}}\bh' \cdot  \nabla^{ \bot} (\Phi^{ \alpha }\gamma'_\rho  (\Psi )\zeta)
\cr
&\qquad \qquad +
\intl_{\R^{2}} \Phi ^\alpha | \nabla \Psi |^2  \gamma_\rho  ''(\Psi ) \zeta + 
\intl_{\R^{2}}   \Phi ^\alpha  \gamma_\rho'(\Psi ) \nabla \Psi\cdot \nabla \zeta. 
\label{8.5}
\end{align}Combining \eqref{8.2} and \eqref{8.5}, we are led  to 
\begin{align}
&-\intl_{\R^{2}}  \Phi ^{ \alpha }\gamma_\rho  (\Psi ) \nabla ^{ \bot} \Psi \cdot  \nabla \zeta
- 
\alpha \Phi^{ \alpha -1} \gamma_\rho  (\Psi ) h^3\zeta
\cr
&\qquad = - \alpha(\alpha -1)\intl_{\R^{2}} | \nabla \Phi|^2 \Phi^{ \alpha -2} \gamma (\Psi )\zeta
+\intl_{\R^{2}}\Phi ^{ \alpha }\gamma_\rho ' (\Psi)  \Delta  \Phi \nabla \Phi \cdot  \nabla ^{ \bot}\zeta
\cr
&\qquad \qquad + \intl_{\R^{2}} \Phi ^\alpha ((\Delta  \Phi)^2 - h^3\Delta \Phi  ) \gamma''_\rho (\Psi ))\zeta
-\intl_{\R^{2}}\bh' \cdot  \nabla^{ \bot} (\Phi^{ \alpha }\gamma'_\rho  (\Psi )\zeta).  
\cr
&\qquad \qquad +
\intl_{\R^{2}} \Phi ^\alpha | \nabla \Psi |^2  \gamma_\rho  ''(\Psi ) \zeta+ 
\intl_{\R^{2}}\Phi ^\alpha     \gamma_\rho'(\Psi ) \nabla \Psi\cdot \nabla \zeta. 
\label{8.6}
\end{align}
 In \eqref{8.6}, first letting $ \alpha \rightarrow 0$, and afterwards letting $ \rho \rightarrow +\infty$, we conclude 
\begin{align}
&- \frac{1}{2}\intl_{\R^{2}} \Psi^2  \nabla ^{ \bot} \Psi \cdot  \nabla \zeta
\cr
&\qquad = -\intl_{\R^{2}}\Psi  \Delta  \Phi \nabla^{ \bot} \Phi \cdot  \nabla \zeta
 + \intl_{\R^{2}}  ((\Delta  \Phi)^2 - h^3\Delta \Phi  ) )\zeta
-\intl_{\R^{2}}\bh' \cdot  \nabla^{ \bot} (\Psi \zeta).  
\cr
&\qquad \qquad +
\intl_{\R^{2}}  | \nabla \Psi |^2   \zeta + 
\intl_{\R^{2}}    \Psi \nabla \Psi\cdot \nabla \zeta. 
\label{8.7}
\end{align}   
Noting that the integral on the left-hand side of \eqref{8.7} vanishes, we deduce the local  energy identity \eqref{8.1}  
from \eqref{8.7} for $ c=0$.  Since in the discussion above  $ \Psi $ can be replaced by $ \Psi -c$ for any $ c\in \R$,  we 
get the assertion of the lemma.  
\hfill \Beweisende   

\vspace{0.3cm}
\hspace{0.5cm}
The following lemma we have used in the  proof of Lemma\,\ref{lem8.1}. 

\begin{lem}
\label{lem8.2}
Let $ \phi \in {\hat W} ^{2,\, 2}(\R^{2})$, and $ \psi \in {\hat W} ^{1,\, 2}(\R^{2})\cap L^\infty(\R^{2})$. Then 
\begin{equation}
\nabla \phi \cdot \nabla^{ \bot} \psi_\var  -(\nabla \phi \cdot \nabla^{ \bot} \psi)_\var    \rightarrow 0  \quad  \text{{\it weakly  in}}\quad  
L^2(\R^{2})\quad  \text{{\it as}}\quad  \var  \rightarrow +\infty.
\label{8.10}
\end{equation}
\end{lem}

{\bf Proof}: By using the absolutely continuity of the Lebesgue measure we see that 
\begin{equation}
\nabla \phi \cdot \nabla^{ \bot} \psi_\var  -(\nabla \phi \cdot \nabla^{ \bot} \psi)_\var    \rightarrow 0  \quad  
\text{{\it a.e. in }}\,\,\,  \R^{2}\quad  \text{{\it as}}\quad  \var  \rightarrow +\infty.  
\label{8.11}
\end{equation} 
Thus, it suffices to show that the $ L^2$-norm of $ \nabla \phi \cdot \nabla^{ \bot} \psi_\var  -(\nabla \phi \cdot \nabla^{ \bot} \psi)_\var$ 
is bounded independently on $ \var $.  To see this, we first calculate for almost everywhere $ x\in \R^{2}$
\begin{align*}
&\nabla \phi (x)\cdot \nabla^{ \bot} \psi_\var(x)  -(\nabla \phi \cdot \nabla^{ \bot} \psi)_\var (x)
\\
&\qquad  = \intl_{B_\var } (\nabla \phi (x)    -\nabla \phi (x-y) )\cdot \nabla^{ \bot} \psi(x-y) \eta_\var(y)  dy 
\\
&\qquad  = -\intl_{B_\var }  \psi(x-y)\intl_{0}^{1} \partial _i \nabla \phi (x- ty) dt y_i\cdot  \nabla ^{ \bot}\eta_\var(y)  dy. 
\end{align*}
Noting that $ | y| \,| \nabla ^{ \bot}\eta_\var(y)|  \lesssim  \var^{ -2} $, and $ \psi \in L^\infty(\R^{2})$, along with Jensen's inequality  we find 
\[
( \nabla \phi (x)\cdot \nabla^{ \bot} \psi_\var(x)  -(\nabla \phi \cdot \nabla^{ \bot} \psi)_\var (x))^2
 \lesssim  \var ^{ -2} \| \psi \|_{ L^\infty}^2 \intl_{B_\var }  \intl_{0}^{1}  | \nabla^2 \phi (x- ty)|^2 dt dy.
\]
Integrating this inequality over $ \R^{2}$, and employing Fubini's theorem,  we obtain 
\begin{align*}
\intl_{\R^{2}} ( \nabla \phi \cdot \nabla^{ \bot} \psi_\var  -(\nabla \phi \cdot \nabla^{ \bot} \psi)_\var )^2
 &\lesssim \var ^{ -2} \| \psi \|_{ L^\infty}^2 \intl_{B_\var }  \intl_{0}^{1} \intl_{\R^{2}} | \nabla^2 \phi (x- ty)|^2 dx dt dy
\\
& =\| \psi \|^2_{ L^\infty}  \| \nabla^2 \phi \|^2_{ L^2}. 
\end{align*}
This completes the proof of the Lemma\,\ref{lem8.2}.  \hfill \Beweisende 

\section{Proof of Theorem \ref{thm1.2}}
\label{sec:-5}
\setcounter{secnum}{\value{section} \setcounter{equation}{0}
\renewcommand{\theequation}{\mbox{\arabic{secnum}.\arabic{equation}}}}

We now consider  the system \eqref{1.2b}, \eqref{1.2a}  in $ \R^{2}$ with general right-hand side. 
From Lemma\,\ref{lem8.1} we infer that 
\begin{align}
&(\Delta \Phi )^2 + | \nabla \Psi |^2 
\cr
&= \nabla \cdot \Big((\Psi - c) \nabla \Psi - (\Psi - c)\Delta \Phi \nabla ^{ \bot} \Phi \Big)
+  h^3 \Delta \Phi - (\Psi -c)(\partial _1 h^1+ \partial _2 h^2)
\label{5.2c}
\end{align}
in $ \R^{2}$  in the sense of distributions. 

\vspace{0.2cm}
\hspace{0.5cm}
Our aim is to prove the following local  regularity result 

\begin{thm}
\label{thm5.1}
Let $ \bh \in \mathcal{M}_{ \rm loc}^{ 2, \mu }$ for some $ \mu >0$. Let $ (\Phi , \Psi )\in 
{\hat W} ^{2,\, 2}(\R^{2})\times {\hat W} ^{1,\, 2}(\R^{2})$ be a strong-weak solution to  \eqref{1.2a}, \eqref{1.2b}. 
Then $ \Psi, \partial _i \Phi \in C^\alpha(\R^{2}) $, $ i=1,2$. 
\end{thm}

\hspace{0.5cm}
The proof of Theorem\,\ref{thm5.1}  is based on Caccioppoli-type 
inequalities as well as a crucial  logarithmic decay estimate.  In what follows we make use of the following notion of a suitable cut off function

\begin{defin}
Given  balls $B_\rho \subset  B_R=B_R(x_0)$, $ 0< \rho < R$,  a function $ \zeta \in C^{\infty}_{\rm c}(B_R)$ is said to be {\it a suitable cut off function for this balls} if $0 \le  \zeta \le 1$ in $ B_R$, $ \zeta \equiv 1$ on $ B_\rho $, and 
\[
| \nabla ^2 \zeta| + | \nabla \zeta |^2  \lesssim  (R-\rho )^{ -2}.    
\]  
\end{defin}

\hspace{0.5cm}
In what follows, let 
 $ \bh \in \bL^2(\R^{2}) $, and let  $ (\Psi, \Phi  )\in  {\hat W} ^{1,\, 2}(\R^{2})\times 
{\hat W} ^{2,\, 2}(\R^{2})$ be a weak-strong solution to  \eqref{1.2b}, \eqref{1.2a}.  Furthermore, 
for a measurable set $ A \subset \R^{2}$ with $ \mes A >0$ we write
\[
(f)_{ A} =  \intmw_{A} f = \frac{1}{\mes A} \intl_{A} f, \quad  f\in L^1(A).  
\]
\ \\
\hspace{0.5cm}
The following lemmas are  an immediate consequence of the local energy 
identity \eqref{8.1}.  
\ \\

{\bf Proof of Theorem\,\ref{thm5.1}}:  
Let $ \zeta \in C^{\infty}_{\rm c}(B_{ R})$ be a cut off function suitable for the balls $ B_{ R/2} 
\subset B_R $. 
In  \eqref{8.1}  we replace $ \zeta $ by $ \zeta ^2$ and set $ c= \Psi _{ B_{ R}  \setminus B_{ R/2}}$. 
This yields  
\begin{align}
& \intl_{B_{ R}} ((\Delta \Phi )^2 + | \nabla \Psi |^2 )\zeta ^2 
\cr
 &\quad  \lesssim R^{ -1} \intl_{B_R  \setminus B_{ R/2}}  | \Psi - (\Psi )_{ B_R  \setminus B_{ R/2}}|\, | \nabla \Psi | \zeta 
\cr
 &\qquad \qquad + R^{ -1}\intl_{B_{ R} \setminus B_{ R/2}}  | \Psi - (\Psi )_{ B_R  \setminus B_{ R/2}}|\, | \Delta  \Phi |
| \nabla \Phi - (\nabla \Phi )_{ B_R  \setminus B_{ R/2}}|   \zeta  
\cr
 &\qquad \qquad + R^{ -1} | (\nabla \Phi) _{ B_{ R}  \setminus B_{ R/2}}|
 \intl_{B_{ R} \setminus B_{ R/2}}  | \Psi - (\Psi )_{ B_R  \setminus B_{ R/2}}|\, | \Delta  \Phi |
\zeta
\cr
 &\qquad \qquad + R^{ -1}\intl_{B_{ R} \setminus B_{ R/2}}  | \Psi - (\Psi )_{ B_R  \setminus B_{ R/2}}|\, | \bh' | 
 \zeta 
\cr
 &\qquad \qquad + \intl_{B_{ R}}  | h^3| \, | \Delta \Phi | \zeta ^2+ | \bh' | \, | \nabla \Psi | \zeta ^2. 
\label{5.0c1}
\end{align}
Then by the aid of H\"older's  inequality, Young's inequality,  and Sobolev-Poincar\'e inequality, we deduce  from 
\eqref{5.0c1}  
\begin{align}
& \intl_{B_{ R}} ((\Delta \Phi )^2 + | \nabla \Psi |^2 )\zeta ^2 
\cr
 &\quad  \lesssim \intl_{B_R  \setminus B_{ R/2}}  | \nabla \Psi |^2
 + \bigg(\intl_{B_{ R} \setminus B_{ R/2}}   | \nabla ^2 \Phi |^2 \bigg)^2
 + | (\nabla \Phi) _{ B_{ R}}|
 \intl_{B_{ R} \setminus B_{ R/2}}  | \nabla \Psi|^2 + | \Delta  \Phi |^2 
 \cr
 &\qquad \qquad + \intl_{B_R} | \bh |^2
\label{5.0c5}.
\end{align}

\hspace{0.5cm}
Next, we provide the estimate of $ \nabla ^2 \Phi $ in term of $ \Delta \Phi $. 
In fact, using integration by parts, we easily get 
\begin{align*}
\intl_{B_R} | \nabla ^2 \Phi |^2 \zeta ^2 &= \intl_{B_R} \partial _i \partial _j\Phi (\partial _i \partial _j\Phi) \zeta ^2 
\\
& = \intl_{B_R} | \Delta \Phi |^2 \zeta ^2   -2\intl_{B_R}  (\partial _j\Phi - (\partial _j\Phi)_{ B_R  \setminus B_{ R/2}})(\partial _i \partial _j\Phi) \zeta \partial _i\zeta
\\
&\qquad \qquad +2\intl_{B_R}  (\partial _j\Phi -  (\partial _j\Phi)_{ B_R  \setminus B_{ R/2}})(\partial _i \partial _j\Phi) \Delta \Phi  \zeta \partial _j\zeta.
\end{align*}
Applying Cauchy-Schwarz's inequality along with Young's inequality and Poincar\'e's inequality, we find 
\begin{align}
\intl_{B_{ R/2}} | \nabla ^2 \Phi |^2  
&  \lesssim  \intl_{B_R} | \Delta \Phi |^2 \zeta ^2   + R^{ -2}\intl_{B_R  \setminus B_{ R/2}}  
| \nabla \Phi - (\nabla \Phi)_{ B_R  \setminus B_{ R/2}})|^2
\cr
 &\lesssim  \intl_{B_R} | \Delta \Phi |^2 \zeta ^2   + \intl_{B_R  \setminus B_{ R/2}}  
| \nabla^2 \Phi |^2.
\label{5.0c4}
\end{align} 
Thus, estimating the first term on the right-hand side of \eqref{5.0c4}  by means of \eqref{5.0c5}, we 
are led to 
\begin{align}
& \intl_{B_{ R/2}} (| \nabla^2  \Phi |^2 + | \nabla \Psi |^2 )
\cr
 &\quad  \lesssim   \bigg[1+| (\nabla \Phi) _{ B_{ R}}| + \intl_{B_{ R} }   | \nabla ^2 \Phi |^2\bigg]
 \intl_{B_{ R} \setminus B_{ R/2}}  | \nabla \Psi|^2 + | \nabla ^2 \Phi |^2 + \intl_{B_R} | \bh |^2\nonumber\\
 & \quad  \lesssim   \bigg[1+| (\nabla \Phi) _{ B_{ R}}| \bigg]
 \intl_{B_{ R} \setminus B_{ R/2}} \left( | \nabla \Psi|^2 + | \nabla ^2 \Phi |^2\right) + R^{\mu}\label{5.0c7}
\end{align}
Thus, by means of \eqref{5.0c7} we are in a position to apply Theorem\,\ref{thm7.1} with $ f= (\partial _1 \Phi , \partial _2\Phi , \Psi )$. Accordingly, 
$ \partial _1\Phi , \partial _2\Phi $ and $\Psi $ are H\"older continuous.  
\hfill \Beweisende 

\vspace{0.5cm}  
{\bf Proof of Theorem\,\ref{thm1.2}}: 
Recalling that $ \bB = (-\partial_2 \Phi , \partial _1, \Psi )$ we obtain the H\"older continuity of $ \bB $. Arguing as in \cite{cha6}, we get the smoothness of $ (\bv, \bB) $.  \hfill \Beweisende

%
%
\section{Proof of  Theorem \ref{thm1.4} } 
\setcounter{secnum}{\value {section}
\setcounter{equation}{0}
\renewcommand{\theequation}{\mbox{\arabic{secnum}.\arabic{equation}}}}


\hspace{0.5cm}
Thanks to Theorem\,\ref{thm1.2} we already know 
that $ \bB $ is smooth, and therefore the following energy identity holds true for all 
 $ \zeta \in C^{\infty}_{\rm c}(\R^{2}) $ and $ \bLambda \in \R^{2} $
\begin{equation}
\intl_{ \R^{2}} | \nabla \bB|^2 \zeta   = \frac{1}{2}\intl_{ \R^{2}} |\bB - \bLambda  |^2 \Delta  \zeta 
+
\intl_{ \R^{2}} (\nabla \times \bB )\times \bB \cdot  (\bB- \bLambda ) \times \nabla\zeta .  
\label{2.0b}
\end{equation}

\hspace{0.5cm}
We define
\[
\mu  (r) := r^{ -2}\intl_{B_r  \setminus B_{ r/2}} | \bB | dx,\quad \quad r >0. 
\] 
By using change of coordinates $ x= r y$, we see that 
\[
\mu (r) = \intl_{B_1  \setminus B_{ 1/2}} | \bB(ry) | dy,\quad \qquad r >0.
\]
By a straightforward arguments we easily get for all $ r> 1$, 
\begin{align*}
\mu ' (r) &= \intl_{B_1  \setminus B_{ 1/2}}  y_j\frac{\partial _j \bB(ry)\cdot \bB (ry) }{| \bB (ry)|}  dy 
\le r^{ -2} \intl_{B_r  \setminus B_{ r/2}} | \nabla \bB | dx
\\
&  \le \,\sqrt{\pi } r^{ -1} \bigg(\intl_{B_r  \setminus B_{ r/2}} | \nabla \bB |^2 dx\bigg) ^{ 1/2} 
  \le \,\sqrt{\pi }  (\log r)' \bigg(\intl_{ \R^{2}  \setminus B_{1/2}} | \nabla \bB |^2 dx\bigg) ^{ 1/2}.
\end{align*}
Thus, the function 
\[
r \mapsto \mu  (r) - \,\sqrt{\pi } \log r \bigg(\intl_{ \R^{2}  \setminus B_{1/2}} | \nabla \bB |^2 dx\bigg) ^{ 1/2}
\]
is non increasing on $ [1, +\infty)$, which implies for all $ r \ge e$
\begin{align}
\mu  (r) & \le \mu (1) +  \,\sqrt{\pi }\log r \bigg(\intl_{ \R^{2}  \setminus B_{1/2}} | \nabla \bB |^2 dx\bigg) ^{ 1/2} 
 \le  C_0 \log(r), 
\label{2.1}
\end{align}
where $ C_0 = \mu (1) +\,\sqrt{\pi }\bigg(\intl_{ \R^{2}  \setminus B_{1/2}} | \nabla \bB |^2 dx\bigg) ^{ 1/2}$. 

\hspace{0.5cm}
Let $ \zeta $ be a cut-off function for $ B_{ r}$ and $ B_{ r/2}$. In \eqref{2.0b} we replace  $ \zeta$ by $\zeta ^2$,  
take  $ \bLambda  =\bB_{ B_r  \setminus B_{ r/2}}$, and integrate by parts.  This gives  
\begin{align*}
\intl_{B_r} | \nabla \bB |^2 \zeta ^2 &=- 2   \intl_{B_r  \setminus B_{ r/2}} \nabla \bB: (\bB- \bB_{ B_r  \setminus B_{ r/2}}) \otimes   \zeta \nabla \zeta 
\\
&\qquad + \intl_{B_r  \setminus B_{ r/2}} (\nabla \times \bB) \times \bB \cdot   (\bB- \bB_{ B_r  \setminus B_{ r/2}}) \times \zeta \nabla \zeta 
\\
& = I_1+I_2. 
\end{align*}
In order to estimate the first integral we use Cauchy-Schwarz inequality together with Poincar\'e's inequality. This gives $ I_1= o(1)$.  
For the estimation of the  second integral $ I_2$  we first write 
\begin{align*}
I_2 &= \intl_{B_r  \setminus B_{ r/2}} (\nabla \times \bB) \times (\bB- \bB_{ B_r  \setminus B_{ r/2}}) \cdot   (\bB- \bB_{ B_r  \setminus B_{ r/2}}) \times \zeta \nabla \zeta 
\\
&\qquad + \intl_{B_r  \setminus B_{ r/2}} (\nabla \times \bB) \times \bB_{ B_r  \setminus B_{ r/2}} \cdot   (\bB- \bB_{ B_r  \setminus B_{ r/2}}) \times \zeta \nabla \zeta 
\\
& = I_{ 21}+ I_{ 22}. 
\end{align*}
Then applying  Cauchy-Schwarz inequality together with Sobolev-Poincar\'e's inequality, we get  $ I_{ 21}= o(1)$ as $r\to +\infty$.  Thus, it only remains to estimate 
$ I_{ 22}$.  By the aid of Cauchy-Schwarz inequality and Poincar\'e's inequality along with \eqref{2.1} we infer 
\begin{equation}
I_{ 22}  \lesssim   \mu (r) \intl_{B_r  \setminus B_{ r/2}} | \nabla \bB |^2 \le C_0 \log(r)\intl_{B_r  \setminus B_{ r/2}} | \nabla \bB |^2.    
\label{2.2}
\end{equation}
Let $ \var >0$, and $ R_0 >0$ be choosen sufficiently large which will be specified below. We now set   $ r = 2^{ k}$, $k\in \N$. 
Let $ \var >0$ be arbitrarily chosen. 
As $ | \nabla \bB|^2 $ is integrable, for every $ m \in \N$ there exists $ k\in  \N, k \ge m$ such that 
\begin{equation}
\intl_{B_{ 2^{k}}  \setminus B_{ 2^{k-1}}} | \nabla \bB |^2 \le  \frac{\var }{k}. 
\label{aaa}
\end{equation}
Otherwise, there exists $m\in \Bbb N$ such that  the reverse inequality of (\ref{aaa}) holds for all $k\geq m$, which leads to the following contradiction
$$
+\infty >\int_{\Bbb R^2} |\nabla \bB|^2 \geq \sum_{k\geq m} \intl_{B_{ 2^{k}}  \setminus B_{ 2^{k-1}}} | \nabla \bB |^2 
\geq \sum_{k\geq m} \frac{\var}{k} =+\infty.
$$
Thus \eqref{2.2} with $ r=2^k$ reads 
\[
I_{ 22}  \lesssim  C_0 \var.  
\] 
As $ \var >0$ can be chosen arbitrarily small,  we conclude that $ \intl_{ \R^{2}} | \nabla \bB |^2=0$ and therefore $ \bB = \const$. 
 \hfill \Beweisende 

\vspace{0.5cm}  
\mbox{\bf Acknowledgements}
Chae was partially supported by NRF grants 2016R1A2B3011647, while Wolf has been supported  by the German Research Foundation (DFG) through the project WO1988/1-1; 612414.

\appendix

\section{Auxiliary Lemmas}
\label{sec:-6}
\setcounter{secnum}{\value{section} \setcounter{equation}{0}
\renewcommand{\theequation}{A.\arabic{equation}}}

\hspace{0.5cm}
By $ {\dot W} ^{m,\, s}(\R^{2})$, $ 1 \le s <+\infty, m\in \N$, we denote the homogeneous Sobolev space of all $ f\in W^{m,\, s}_{ \rm loc}(\R^{2})$ with $ D^\alpha  f_{ B_1} = \frac{1}{\mes B_1} \intl_{ \R^{2}} D^\alpha  f=0 $ for all 
$ | \alpha | \le m-1$ such that $ D^\alpha  f \in  L^s(\R^{2})$.

\begin{lem}
Let $ B_1, B_2 \in {\dot W} ^{1,\, 2}(\R^{2})$ with $ \partial _1 B^1+ \partial _2 B^2 =0$ in $ \R^{2}$. Then there exists $  \Phi \in {\dot W} ^{2,\, 2}(\R^{2}) $, 
such that $ (B^1, B^2)^{ \top} = \nabla ^{ \bot} \Phi $, where $\nabla ^{ \bot} \Phi = (\partial _2 \Phi , - \partial _1 \Phi )^{ \top} $.  
\end{lem}

{\bf Proof}: We consider the equation 
\begin{equation}
-\Delta \Phi = \partial _1 B^2 - \partial _2 B^1\quad  \text{ in}\quad  \R^{2}. 
\label{3.2}
\end{equation}
Let $\Phi \in  {\dot W} ^{2,\, 2}(\R^{2}) $ denote the unique weak solution   to \eqref{3.2}.  
Set $ A= \nabla ^{ \bot} \Phi $. Then $ \partial _1 A^2 - \partial _2 A^1= -\Delta \Phi = \partial _1 B^2 - \partial _2 B^1$. 
Thus, there exists $ p \in  W_{ \rm loc} ^{1,\, 2}(\R^{2}) $, such that $ (B^1-A^1, B^2-A^2) = (\partial _1 p, \partial _2 p )$. 
Since $ \partial _1 (B^1-A^1) + \partial _2 (A^2-B^2)$ it follows $ \Delta p =0$.  Noting that $ \nabla p$ has logarithmic growth at 
infinity we get $ \nabla p = \const$.  Eventually, replacing $ \Phi $ by $ \Phi + Q$, where $ Q$ is a polynomial of degree $ \le 1$ we may assume that $ (B^i- A^i)_{ B_1} = 0$ \, $ (i=1,2)$, which gives $ (B^1, B^2) = (A^1, A^2) =\nabla ^{ \bot} \Phi $. 
 \hfill \Beweisende \\

\begin{lem}
\label{lemA.7} 
Let $ \phi: [0, 1) \rightarrow  [0, +\infty)$ be a non decreasing function. Assume there exists $ c= \const >0$ such that for all $ 0<r< \frac{1}{2}$
\begin{equation}
\phi (r/2) \le \frac{c[\log r^{ -1}]^{\frac12} }{1+ c[\log r^{ -1}]^{\frac12} }\phi (r) +  c r^{\mu}. 
\label{6.27}
\end{equation}
Then for every $ \alpha  >1$ there holds for all $ 0<r< \frac{1}{2}$ 
\begin{equation}
\phi (r)  \lesssim   \frac{1}{[\log r^{ -1}]^{\alpha }},
\label{6.28}
\end{equation}
where the hidden constant in \eqref{6.28} depend only on $ c, q, \alpha $ and $ \mu $.  

\end{lem}

{\bf Proof}: Clearly \eqref{6.27} with $ r= 2^{ -k}$  reads 
\begin{equation}
\phi (2^{ -k}) \le \frac{c k^\frac12 }{1+ c  k^\frac12 }\phi (2^{ -k+1}) +  c 2^{ -k\mu +\mu }. 
\label{6.29}
\end{equation}
Let $ n\in \N, n \ge 2$. Iterating \eqref{6.29} from $ k=n^2 $  to $ k= n+1$, we obtain 
\begin{align}
\phi (2^{-n^2})& \le  \Big(\prod_{ k=n+1}^{ n^2} \frac{c k^\frac12  }{1+ c k^\frac12  }\Big)\phi (2^{ n}) +  \frac{c 2^{ -n \mu }}{1- 2^{-\mu}}.
\label{6.30}
\end{align}
Furthermore, estimating 
\begin{align*}
\log \prod_{ k=n+1}^{ n^2} \frac{c k^\frac12  }{1+ c k^\frac12  }
&=\sum_{k=n+1}^{n^2} \log \Big(1- \frac{1}{1+ ck ^\frac12}\Big)
\le  - \sum_{k=n+1}^{n^2} \frac{1}{1+ ck^\frac12  } \lesssim  - n^{ \frac12},
\end{align*}
we find that 
\begin{equation}
\prod_{ k=n+1}^{ n^2} \frac{c k^\frac12 }{1+ c k^\frac12 }  \lesssim  e^{- n^{ \frac12 } }.
\label{6.31}
\end{equation}
Noting that $  \frac{ 2^{ -n \mu }}{1- 2^{-\mu}}  \lesssim e^{- n^{ \frac12 } } $, we deduce from \eqref{6.30} along with \eqref{6.31} that 
 \begin{equation}
 \phi (2^{-n^2})  \lesssim  e^{- n^{ \frac12 } }
 \label{2.26}
  \end{equation}
  Let $ \alpha >1$ be arbitrarily chosen.  Clearly,  there exists $ n_0 \in \N$ such that for all $ n \ge n_0$
\[
2\alpha  \log n \le n^{ \frac12 }\quad  \Longrightarrow\quad  e^{- n^{ \frac12} } \le \frac{1}{[\log 2^{ n^2}]^\alpha }. 
\]  
Let $ 0< r< \frac{1}{2}$. Then there exists  unique $ n\in \N$ such that $2^{ -(n+1)^2} < r \le 2^{ -n^2} $. 
In particular, $ (n+1)^2 > \frac{\log r^{ -1}}{\log 2}$, and \eqref{2.26} yields 
\[
\phi (r)  \le \phi (2^{ -(n+1)^2})     \lesssim  
 \frac{1}{[\log r^{ -1}]^{\alpha }}.
\] 
Whence, the claim.  \hfill \Beweisende \\

\begin{lem}
\label{lemA.5}
Let $ f\in {\hat W} ^{1,\, 2}(\R^{2})$. Suppose there exists $ \alpha >1$ and a constant $ c_0>0$ such that for all $ 0< r < \frac{1}{2}$ 
\begin{equation}
\bigg(\intl_{B_r} | \nabla f|^2 \bigg)^{ \frac{1}{2}} \le  \frac{c_0}{[\log r^{ -1}]^\alpha }.  
\label{6.17}
\end{equation}
Then 
\begin{equation}
\supl _{ 0< r< 1} | (f)_{ B_r}| <+\infty .
\label{6.18}\end{equation}
\end{lem}

{\bf Proof}:  We set $ {\D \mu (r) = \intmw_{B_r} f} $.  We first claim the following inequality for  all $0<r<R< +\infty$.
\begin{equation}
\label{6.18a}
| \mu (r)| \le | \mu (R)|   +\log \frac{R}{r}  \left( \pi ^{-1} \int_{ B_R} |\nabla f|^2 \right) ^{1/2}.
\end{equation}
Indeed, since $ \mu(r)= \intmw_{\,B_1} f(ry)dy$,  we have
$$
\mu '(r)=\intmw_{B_1} y\cdot \nabla f(ry) dy = \intmw_{B_r} \frac{x}{r}\cdot \nabla f dx \geq -\intmw_{B_r} |\nabla f|
\geq -\frac{1}{r}\left( \pi^{-1} \int_{B_R} |\nabla f|^2\right)^{1/2}.
$$
This implies that
$$
\frac{d}{dr} \left[ \mu(r)+\log r \left( \pi ^{-1} \int_{B_R} |\nabla f|^2 \right)^{1/2} \right] \geq 0,
$$
which after integration over $[r, R]$ provides us with
\begin{equation}
\label{6.18b}
\mu (r)\leq \mu (R) +\log \frac{R}{r}  \left( \pi ^{-1} \int_{ B_R} |\nabla f|^2 \right) ^{1/2}.
\end{equation}
Similarly, we have
$$
\mu '(r)\leq  \frac{1}{r}\left( \pi^{-1} \int_{B_R} |\nabla f|^2\right)^{1/2},
$$
which leads to the opposite inequality,
\begin{equation}
\label{6.18c}
\mu (r)\geq \mu (R) -\log \frac{R}{r}  \left( \pi ^{-1} \int_{ B_R} |\nabla f|^2 \right) ^{1/2}.
\end{equation}
Combining (\ref{6.18b}) and (\ref{6.18c}), we obtain (\ref{6.18a}) as claimed. 

\hspace{0.5cm}
For $ k \in \N$, $ k \ge 2$   we insert $ r= 2^{ -k-1}$ and $ R= 2^{ -k}$ in  \eqref{6.18a}. This yields 
\begin{equation}
| \mu (2^{ -k-1})| \le | \mu ( 2^{ -k})|   + \log 2 \frac{c_0 \pi^{- \frac{1}{2}}}{ k^\alpha [\log 2]^\alpha }
\le | \mu ( 2^{ -k})|   +  \frac{c_0 }{ k^\alpha}. 
\label{6.20}
\end{equation}
Iterating \eqref{6.20} $ n-1$-times from $ k=n$ to $ k=2$,  we arrive at 
\[
| \mu (2^{ -n-1})| \le | \mu (2^{ -2})| + c_0 \sum_{k=2}^{n} k^{ -\alpha }  \lesssim  
\zeta (\alpha ) .
\]  
Here $ \zeta $  stands for Riemann's  Zeta-function, and $ \zeta (\alpha ) = \sum_{k=1}^{\infty} k^{ -\alpha }$.

\hspace{0.5cm}
Now, let $ 0< r < \frac{1}{2}$ be arbitrarily chosen. There exists  unique $ n \in \N$ such that $ 2^{ -n-2} < r \le 2^{ -n-1}$. 
In particular, $ 2n \ge n +2\ge \frac{\log r^{ -1}}{\log 2} \ge \log r^{ -1}$. Thus, the  inequality above yields 
$
| \mu (r)|  \lesssim 
\zeta (\alpha ).
$
This completes the proof of \eqref{6.18}.  \hfill \Beweisende  \\

\begin{lem}
\label{lemA.8} 
Let $ \phi: [0, 1) \rightarrow  [0, +\infty)$ be a non decreasing function. Assume there exists a constants  $\theta, \mu  \in (0,1)$ such that for all $ 0<r< \frac{1}{2}$
\begin{equation}
\phi (r/2) \le \theta \phi (r) +  c r^{\mu}. 
\label{6.32}
\end{equation}
Then for every $ \alpha  >1$ there holds for all $ 0<r< \frac{1}{2}$ 
\begin{equation}
\begin{cases}
\phi (r)  \lesssim   r^{ \alpha }\quad  &\text{ if}\quad   \theta \not= 2^{ -\mu },\quad  \alpha = 
\min \Big\{-\frac{\log \theta }{\log 2}, 2\Big\}, 
\\[0.3cm]
\phi (r)  \lesssim  (\log r^{ -1})  r^{ \mu }\quad  &\text{ if}\quad  \theta = 2^{ -\mu }
\end{cases}
\label{6.33}
\end{equation}
where the hidden constant in \eqref{6.28} depends only on $ c, \theta $ and $ \mu $.  

\end{lem}

{\bf Proof}:  The inequality \eqref{6.32} with $ r= 2^{ -k}$ reads 
\begin{equation}
\phi (2^{ -k-1}) \le \theta \phi (2^{ -k}) + c 2^{ -k\mu }. 
\label{6.34}
\end{equation}  
Iterating \eqref{6.34} $ n$-times,  we obtain 
\begin{equation}
\phi (2^{ -n}) \le c 2^{ -n\mu } + c \theta 2^{ -(n-1)\mu } +c \theta^2 2^{ -(n-2)\mu } + \ldots + c \theta^{ n-1} 
2^{\mu }+c \theta^{ n}.
\label{6.35}
\end{equation} 
In case $ 2^{ -\mu } < \theta $ we deduce from \eqref{6.35}
\begin{align}
\phi (2^{-n}) &\le c \theta^{ n} \bigg[ \Big(\frac{2^{ -\mu }}{\theta }\Big)^n+ \Big(\frac{2^{ -\mu }}{\theta }\Big)^{ n-1}+ 
\ldots + 1\bigg]
\le   \frac{c\theta^n }{\theta - 2^{ -\mu }}.
\label{6.36}
\end{align} 
On the contrary, in case  $ 2^{ -\mu } > \theta $ \eqref{6.35} yields 
\begin{align}
\phi (2^{-n}) &\le c 2^{-\mu  n} \bigg[ 1+ \frac{\theta }{2^\mu }+ \ldots + 
\Big(\frac{\theta }{2^{ -\mu }}\Big)^{ n-1}+ \Big(\frac{\theta }{2^{ -\mu }}\Big)^{ n}\bigg]
\le   \frac{c 2^{ -\mu n } }{  2^{ -\mu }- \theta}.
\label{6.37}
\end{align} 
Thus, from \eqref{6.36} and \eqref{6.37} we deduce that  for all $ n\in \N$ 
\begin{equation}
\phi (2^{-n})  \lesssim   2^{ -n \alpha },\quad  \alpha = \min \Big\{\frac{\log \theta }{\log 2}, \mu \Big\}. 
\label{6.38}
\end{equation}
This implies \eqref{6.33} in case $ \theta \not= 2^{ -\mu }$. 

\hspace{0.5cm}
 In case $ \theta = 2^{ -\mu }$ we infer from \eqref{6.35} 
\begin{equation}
\phi (2^{-n})  \lesssim   n 2^{ -n\mu}. 
\label{6.39}
\end{equation} 
Noting that $ n = \frac{ \log 2^{n}}{\log 2}$, we get \eqref{6.33} in case $ \theta = 2^{ -\mu }$.  \hfill \Beweisende

\bibliographystyle{siam}
\bibliography{ChaeWolf2015b}


\end{document}